\documentclass[10pt]{article}

\usepackage{amsmath}%
\usepackage{amsfonts}%
\usepackage{amssymb}%
\usepackage{multicol}%
\usepackage{graphicx}%
\usepackage{caption}
\usepackage{float}
\usepackage[titletoc]{appendix}

\newcommand{\hyp}[5]{\,\mbox{}_{#1}F_{#2}\!\left(\genfrac{}{}{0pt}{}{#3}{#4};#5\right)}
\renewcommand\Re{\operatorname{Re}}
\numberwithin{equation}{section}

\usepackage{geometry}
\newtheorem{theorem}{Theorem}

\newtheorem{Remark}[theorem]{Remark}

\begin{document}

\title{Some definite integrals involving Jacobi polynomials.}
\author{Enno Diekema \footnote{email adress: e.diekema@gmail.com}}
\maketitle

\begin{abstract}
\noindent
Szmytkowski derived a certain integral with Gegenbauer polynomials. A natural generalization is to derive lookalike integrals with Jacobi polynomials. Six methods are treated to derive the first integral. The first method should be enough to prove the first integral, but by the other methods there arises remarkable formula such as par example a zero-balanced $F_3$ Appell function which can be converted into a $_2F_1$ hypergeometric function. Another three integrals complete the paper.
\end{abstract}

\section{Introduction}
\setlength{\parindent}{0cm}
During our investigations on fractional calculus \cite {4} we met the integral
\begin{equation}
I=\int_0^1\dfrac{(1-t)^\alpha(1+t)^\beta}{t^\lambda}P_n^{(\alpha,\beta)}(t)dt \qquad \text{with}\qquad  \lambda>0\ \ \text{and}\ \ \Re(\alpha,\beta)>-1
\label{1}
\end{equation}
where $P_n^{(\alpha,\beta)}(t)$ are the Jacobi polynomials. Searching in the literature and on the internet does only found an integral from Szmytkowski with Gegenbauer polynomials. He derived the integral formula \cite[2.7]{1} 
\[
\dfrac{n+a}{a}\int_x^1\dfrac{(1-t^2)^{a-1/2}}{(t-x)^{\kappa+1/2}}C^{(a)}_n(t)dt=\dfrac{\sqrt{\pi}\Gamma(n+2a)(n+a)\Gamma(1/2-\kappa)}{2^{a-1/2}n!\Gamma(a+1)}(1-x^2)^{(a-\kappa)/2}
P^{\kappa-a}_{n+a-1/2}(x)
\]
with conditions: $\Re(a)>-\dfrac{1}{2}$,\ \ $\Re(\kappa)<\dfrac{1}{2}$ and $-1<x<1$. $C^a_n(t)$ are the Gegenbauer polynomials while $P^{\kappa-a}_{n+a-1/2}(x)$ is the associated Legendre function of the first kind $\text{\em on the cut}$.

\

Substitution of $x=0$ gives
\begin{equation}
\dfrac{n+a}{a}\int_0^1\dfrac{(1-t^2)^{a-1/2}}{t^{\kappa+1/2}}C^a_n(t)dt=\dfrac{\sqrt{\pi}\Gamma(n+2a)(n+a)\Gamma(1/2-\kappa)}{2^{a-1/2}n!\Gamma(a+1)}P^{\kappa-a}_{n+a-1/2}(0).
\label{2}
\end{equation}
Supposing $\lambda=\kappa+\dfrac{1}{2}$ gives $\lambda<1$. Then \eqref{2} reduces to
\begin{equation}
\dfrac{n+a}{a}\int_0^1\dfrac{(1-t^2)^{a-1/2}}{t^\lambda}C^a_n(t)dt=\dfrac{\sqrt{\pi}\Gamma(n+2a)(n+a)\Gamma(1-\lambda)}{2^{a-1/2}n!\Gamma(a+1)}P^{\lambda-a -1/2}_{n+a-1/2}(0).
\label{3}
\end{equation}
For the associated Legendre function of the first kind we have \cite[8.756(1)]{2}:
\[
P^{\lambda-a-1/2}_{n+a-1/2}(0)=
\dfrac{2^{\lambda-a-1/2}\sqrt{\pi}}{\Gamma\left(\dfrac{n+2a-\lambda+2}{2}\right)
\Gamma\left(\dfrac{2-\lambda-n}{2}\right)}.
\]
Then we obtain from \eqref{3}
\begin{equation}
\int_0^1\dfrac{(1-t^2)^{a-1/2}}{t^\lambda}C^a_n(t)dt=
\dfrac{\pi\Gamma(n+2a)\Gamma(1-\lambda)}{\Gamma(n+1)\Gamma(a)}\dfrac{2^{\lambda-2a}}{\Gamma\left(\dfrac{n+2a-\lambda+2}{2}\right)	\Gamma\left(\dfrac{2-\lambda-n}{2}\right)}.
\label{4}
\end{equation}
A natural generalization is the integral with the Jacobi polynomials and the matching weight function. The orthogonality bounds for these polynomials are $-1$ and $1$. But because the singularity for $t=0$ we use the bounds $0$ and $1$. Some other integrals can be derived from the basic integral.

In this paper the Beta function is defined as
\[
B(a,b)=\dfrac{\Gamma(a)\Gamma(b)}{\Gamma(a+b)}
\]

\

\section{Overview of the integrals concerning the Jacobi polynomials}
In this paper we treat the following theorems
\begin{theorem}
\begin{equation}
\int_{-1}^1\dfrac{(1-t)^\alpha(1+t)^\beta}{(z-t)^\lambda}P_n^{(\alpha,\beta)}(t)dt=
\dfrac{2^{\alpha+\beta+n+1}(\lambda)_n}{(z-1)^{n+\lambda}n!}B(\alpha+n+1,\beta+n+1)
\hyp21{\alpha+n+1,n+\lambda}{\alpha+\beta+2n+2}{\dfrac{2}{1-z}}
\label{2.000}
\end{equation}
with $|z|>1$. 	
\end{theorem}
This integral is a generalization of \cite[16.4.(4)]{7} where the integral is given for $\lambda=1$.
\begin{theorem}
\begin{equation}
\int_x^1\dfrac{(1-t)^\alpha(1+t)^\beta}{(t-x)^\lambda}P_n^{(\alpha,\beta)}(t)dt=
2^\beta\dfrac{\Gamma(n+\alpha+1)\Gamma(1-\lambda)}{\Gamma(n+1)\Gamma(\alpha-\lambda+2)}
(1-x)^{\alpha+1-\lambda}\hyp21{\alpha+n+1,-\beta-n}{\alpha-\lambda+2}{\dfrac{1-x}{2}}
\label{2.00}
\end{equation}

with $-1<x<1$.
\end{theorem}

\begin{theorem}
\begin{equation}
\int_0^1\dfrac{(1-t)^\alpha(1+t)^\beta}{t^\lambda}P_n^{(\alpha,\beta)}(t)dt
=2^\beta\dfrac{\Gamma(n+\alpha+1)\Gamma(1-\lambda)}{\Gamma(n+1)\Gamma(\alpha-\lambda+2)}
\hyp21{\alpha+n+1,-\beta-n}{\alpha-\lambda+2}{\dfrac{1}{2}}
\label{2.001}
\end{equation}
\end{theorem}
This is a special case of Theorem 2 with $x=0$.

\begin{theorem}
\begin{align}
\int^x_{-1}\dfrac{(1-t)^\alpha(1+t)^\beta}{(t-x)^\lambda}P_n^{(\alpha,\beta)}(t)dt
&=\dfrac{2^{\alpha+\beta+n+1}}{(1-x)^{n+\lambda}}\binom{-\lambda}{n}
B(\alpha+n+1,\beta+n+1)\hyp21{\alpha+n+1,n+\lambda}{\alpha+\beta+2+2n}{\dfrac{2}{1-x}}- \nonumber \\
&-2^\beta\dfrac{\Gamma(n+\alpha+1)\Gamma(1-\lambda)}{\Gamma(n+1)\Gamma(\alpha-\lambda+2)}
(1-x)^{\alpha+1-\lambda}\hyp21{\alpha+n+1,-\beta-n}{\alpha-\lambda+2}{\dfrac{1-x}{2}}
\label{2.003}
\end{align}
with $-1<x<1$. 
\end{theorem}
For all these integrals the conditions are $\Re(\alpha,\beta)>-1$ and $0<\lambda \leq 1$.

\

\section{Proof of the theorems}

\textbf{\fontsize{10.5}{12.5}\selectfont Proof of Theorem 1.}

We treat several methods to compute the integral.
For the first method we use the formula of Rodrigues for the Jacobi polynomials. This is a very direct method. The next methods are given because there arises many remarkable formulas. For the second method we use an integral given in \cite[2.22.4.11]{3}. This integral can be converted into the desired one. The third method uses also a formula of a known integral \cite[(2.22.4.9)]{3}. Because there are no proofs of these formula in the reference the first method can be used as a proof of these formula. The fourth method uses again the formula of Rodrigues for the Jacobi polynomials. The fifth method is a variation of the fourth one. The sixth method is a so-called brute force method. The Jacobi polynomials are written as a summation. Interchanging the integral and the summation gives after a lot of manipulations with the Gamma functions and the Pochhammer symbols the desired result.

\

\textbf{First method.}

This first method starts with the Rodrigues formula for the Jacobi polynomials
\[
P_n^{(\alpha,\beta)}(t)=\dfrac{(-1)^n}{2^n\, n!}\dfrac{1}{(1-t)^\alpha(1+t)^\beta}
\dfrac{d^n}{dt^n}\big((1-t)^\alpha(1+t)^\beta(1-t^2)^n \big).
\]
Application to the integral \eqref{2.000} gives
\[
I=\int_{-1}^1\dfrac{(1-t)^\alpha(1+t)^\beta}{(z-t)^\lambda}P_n^{(\alpha,\beta)}(t)dt=
\dfrac{(-1)^n}{2^n\, n!}\int_{-1}^1\dfrac{1}{(z-t)^\lambda}\dfrac{d^n}{dt^n}
\big((1-t)^{\alpha+n}(1+t)^{\beta+n}\big)dt.
\]
$n$-times partial integration gives
\[
I=\dfrac{2^{\alpha+\beta+n+1}(\lambda)_n}{(z+1){\lambda+n}n!}B(\alpha+n+1,\beta+n+1)
\hyp21{\beta+1+n,\lambda+n}{\alpha+\beta+2n+2}{\dfrac{2}{z+1}}
\]
Using a standard Gaussian transformation gives
\[
\int_{-1}^1\dfrac{(1-t)^\alpha(1+t)^\beta}{(z-t)^\lambda}P_n^{(\alpha,\beta)}(t)dt=
\dfrac{2^{\alpha+\beta+n+1}(\lambda)_n}{(z-1){\lambda+n}n!}B(\alpha+n+1,\beta+n+1)
\hyp21{\alpha+1+n,\lambda+n}{\alpha+\beta+2n+2}{\dfrac{2}{1-z}}
\]
This proves the theorem. $\square$

\

\textbf{Second method.}

Our starting formula is the integral \cite[(2.22.4.11)]{3} with $a=1$,\ $\theta=\lambda$ and $\beta=\alpha+1$. 
\begin{multline*}
\int_{-1}^1\dfrac{(1-t)^\alpha(1+t)^\sigma}{(z-t)^\lambda}P_n^{(\rho,\sigma)}(t)dt= \\
=\dfrac{(\rho-\alpha)_n}{n!}B(\alpha+1,\sigma+n+1)2^{\alpha+1+\sigma}\dfrac{1}{(z-1)^\lambda}
\hyp32{\lambda,\alpha+1,\alpha+1-\rho}{\alpha+1-\rho-n,\alpha+\sigma+n+2}{\dfrac{2}{1-z}}
\end{multline*}
with $\alpha,\sigma>-1$, $\lambda \leq 1$ and $arg(z^2-1)<\pi$. Setting $\sigma=\beta$ and $\rho=\alpha+\epsilon$ with $\epsilon \rightarrow 0$ gives
\begin{multline*}
I=\int_{-1}^1\dfrac{(1-t)^\alpha(1+t)^\beta}{(z-t)^\lambda}P_n^{(\alpha,\beta)}(t)dt= \\ 
=\dfrac{(\epsilon)_n}{n!}B(\alpha+1,\beta+n+1)2^{\alpha+\beta+1}\dfrac{1}{(z-1)^\lambda}
\hyp32{\lambda,\alpha+1,1}{1-\epsilon-n,\alpha+\beta+n+2}{\dfrac{2}{1-z}}
\end{multline*}
Rewriting this equation gives
\begin{multline*}
I=\dfrac{1}{n!}\dfrac{\Gamma(n+\epsilon)\Gamma(1-\epsilon-n)}{\Gamma(\epsilon)}
B(\alpha+1,\beta+n+1)2^{\alpha+\beta+1} \\
\dfrac{1}{(z-1)^\lambda}\dfrac{1}{\Gamma(1-\epsilon-n)}
\hyp32{\lambda,\alpha+1,1}{1-\epsilon-n,\alpha+\beta+n+2}{\dfrac{2}{1-z}}
\end{multline*}
Using 
\[
\lim_{\epsilon \rightarrow 0}
\dfrac{\Gamma(n+\epsilon)\Gamma(1-\epsilon-n)}{\Gamma(\epsilon)}=(-1)^n
\]
gives
\begin{equation}
I=\dfrac{(-1)^n}{n!}
B(\alpha+1,\beta+n+1)2^{\alpha+\beta+1} \\
\dfrac{1}{(z-1)^\lambda}\dfrac{1}{\Gamma(1-\epsilon-n)}
\hyp32{\lambda,\alpha+1,1}{1-\epsilon-n,\alpha+\beta+n+2}{\dfrac{2}{1-z}}
\label{2.001}
\end{equation}
We use the following property \cite[Lemma 2]{6} with $M$ a non-negative integer
\begin{multline}
\dfrac{1}{\Gamma(-M)} \ _{p+1}F_p
\left(\begin{array}{l}
a_0,\dots,a_p \\
-M,b_2,\dots,b_p
\end{array};z\right)= \\
=\dfrac{z^{M+1}(a_0)_{M+1}\dots(a_p)_{M+1}}{\Gamma(M+2)(b_2)_{M+1}\dots(b_p)_{M+1}}
\ _{p+1}F_p
\left(\begin{array}{l}
a_0+M+1,\dots,a_p+M+1 \\
M+2,b_2+M+1,\dots,b_p+M+1
\end{array};z\right).
\label{2.002}
 \end{multline}
Applying this property with $M=n-1+\epsilon$ we get for \eqref{2.001} after some simplification
\[
I=\dfrac{(\lambda)_n}{n!}B(\alpha+n+1,\beta+n+1)2^{\alpha+\beta+1}\dfrac{1}{(z-1)^{n+\lambda}}
\hyp21{n+\lambda,\alpha+n+1}{\alpha+\beta+2n+2}{\dfrac{2}{1-z}}
\]
This proves the theorem. $\square$

\pagebreak
\textbf{Third method.}

Our starting formula is the integral \cite[(2.22.4.9)]{3}
\begin{align}
&\int_{-b}^bt^m(b-t)^\alpha(b+t)^\beta P_n^{(\alpha,\beta)}\left(\dfrac{t}{b}\right)dt=I_m \label{2.1} \\
&\text{with} \nonumber \\
&b>0, \ \ \qquad \Re(\alpha,\beta)>-1 \nonumber \\
&I_m=0, \qquad m=0,1,2,\dots,n-1 \nonumber \\
&I_n=B(\alpha+n+1,\beta+n+1)(2b)^{\alpha+\beta+n+1} \nonumber \\
&I_m=\binom{m}{n}B(\alpha+n+1,\beta+n+1)2^{\alpha+\beta+n+1}
b^{\alpha+\beta+m+1}\hyp21{n-m,\alpha+n+1}{\alpha+\beta+2n+2}{2} \qquad m>n. \nonumber
\end{align}

We use the next formula
\begin{equation}
\dfrac{1}{(z-t)^\lambda}=\dfrac{1}{z^\lambda}\dfrac{1}{\left(1-\dfrac{t}{z}\right)^\lambda}=\dfrac{1}{z^\lambda}
\sum_{k=0}^\infty \dfrac{(\lambda)_k}{k!}\left(\dfrac{t}{z}\right)^k
\label{2.2a}
\end{equation}
with conditions $\left|\dfrac{t}{z}\right|<1$. Because $-1<t<1$ there follows $|z|>|t|$.

\ 

The integral $I$ can be computed using \eqref{2.1} with $b=1$ and \eqref{2.2a}
\begin{align}
I&=\int_{-1}^1\dfrac{(1-t)^\alpha(1+t)^\beta}{(z-t)^\lambda}P_n^{(\alpha,\beta)}(t)dt \nonumber \\
&=\dfrac{1}{z^\lambda}\int_{-1}^1(1-t)^\alpha(1+t)^\beta P_n^{(\alpha,\beta)}(t)
\sum_{k=0}^\infty \dfrac{(\lambda)_k}{k!}\left(\dfrac{t}{z}\right)^kdt \nonumber \\
&=\dfrac{1}{z^\lambda}\sum_{k=0}^\infty \dfrac{(\lambda)_k}{k!z^k}
\int_{-1}^1 t^k(1-t)^\alpha(1+t)^\beta P_n^{(\alpha,\beta)}(t)dt= \nonumber\\
&=\dfrac{2^{\alpha+\beta+n+1}B(\alpha+n+1,\beta+n+1)}{z^\lambda n!}
\sum_{k=0}^\infty \dfrac{(\lambda)_k}{z^k\Gamma(k-n+1)}\hyp21{n-k,\alpha+n+1}{\alpha+\beta+2n+2}{2}.
\label{2.6}
 \end{align}
The interchanging of the integral and the summation is allowed because of the convergence of the summation. For $k \leq n-1$ the integral and so the summation is equal $0$. Then the lowest value of $k$ is $k=n$ and \eqref{2.6} becomes
\[
I=\dfrac{2^{\alpha+\beta+n+1}B(\alpha+n+1,\beta+n+1)}{z^\lambda n!}
\sum_{k=n}^\infty \dfrac{(\lambda)_k}{z^k\Gamma(k-n+1)}\hyp21{n-k,\alpha+n+1}{\alpha+\beta+2n+2}{2}.
\]
Setting $k=m+n$ results in
\begin{equation}
I=\dfrac{2^{\alpha+\beta+n+1}B(\alpha+n+1,\beta+n+1)}{z^{n+\lambda}\Gamma(n+1)\Gamma(\lambda)}
\sum_{m=0}^\infty \dfrac{\Gamma(m+n+\lambda)}{\Gamma(m+1)}\dfrac{1}{z^m}
\hyp21{-m,\alpha+n +1}{\alpha+\beta+2n+2}{2}.
\label{2.12}
\end{equation}
The hypergeometric function can be written as a summation. We get
\[
\hyp21{-m,\alpha+\beta+1}{\alpha+\beta+2n+2}{2}=
\dfrac{\Gamma(\alpha+\beta+2n+2)}{\Gamma(-m)\Gamma(\alpha+n+1)}
\sum_{k=0}^m\dfrac{\Gamma(k-m)\Gamma(k+\alpha+n+1)}{\Gamma(k+\alpha+\beta+2n+2)}
\dfrac{2^k}{k!}.
\]
After substitution in \eqref{2.12} and using the definition of the Beta function we obtain
\[
I=\dfrac{2^{\alpha+\beta+n+1}\Gamma(\beta+n+1)}{z^{n+\lambda}\Gamma(n+1)\Gamma(\lambda)}
\sum_{m=0}^\infty\dfrac{\Gamma(m+n+\lambda)}{\Gamma(m+1)\Gamma(-m)}\dfrac{1}{z^m}
\sum_{k=0}^m\dfrac{\Gamma(k-m)\Gamma(k+\alpha+n+1)}{\Gamma(k+\alpha+\beta+2n+2)}
\dfrac{2^k}{k!}.
\]
Rewriting the double summation gives
\[
I=\dfrac{2^{\alpha+\beta+n+1}\Gamma(\beta+n+1)}{z^{n+\lambda}\Gamma(n+1)\Gamma(\lambda)}
\sum_{k=0}^\infty\dfrac{\Gamma(k+\alpha+n+1)}{\Gamma(k+\alpha+\beta+2n+2)}\dfrac{2^k}{k!}
\sum_{m=k}^\infty\dfrac{\Gamma(k-m)\Gamma(m+n+\lambda)}{\Gamma(m+1)\Gamma(-m)}\dfrac{1}{z^m}.
\]
For the Gamma functions there is the following property
\begin{equation}
\dfrac{\Gamma(k-m)}{\Gamma(-m)}=(-1)^k\dfrac{\Gamma(1+m)}{\Gamma(1+m-k)}
\label{2.12aa}
\end{equation}
Application yields
\[
I=\dfrac{2^{\alpha+\beta+n+1}\Gamma(\beta+n+1)}{z^{n+\lambda}\Gamma(n+1)\Gamma(\lambda)}
\sum_{k=0}^\infty\dfrac{\Gamma(k+\alpha+n+1)}{\Gamma(k+\alpha+\beta+2n+2)}\dfrac{(-2)^k}{k!}
\sum_{m=k}^\infty\dfrac{\Gamma(m+n+\lambda)}{\Gamma(m+1-k)}\dfrac{1}{z^m}.
\]
The last summation is well-known
\[
\sum_{m=k}^\infty\dfrac{\Gamma(m+n+\lambda)}{\Gamma(m+1-k)}\dfrac{1}{z^m}=
\Gamma(k+n+\lambda)\left(\dfrac{z}{z-1}\right)^{\lambda+k+n}\dfrac{1}{z^k}.
\]
Substitution gives
\begin{equation}
I=\dfrac{2^{\alpha+\beta+n+1}\Gamma(\beta+n+1)}{(z-1)^{n+\lambda}\Gamma(n+1)\Gamma(\lambda)}
\sum_{k=0}^\infty\dfrac{\Gamma(k+\alpha+n+1)\Gamma(k+n+\lambda)}{\Gamma(k+\alpha+\beta+2n+2)}\dfrac{1}{k!}
\left(\dfrac{2}{1-z}\right)^k.
\label{2.6a}
\end{equation}
Writing the summation as a hypergeometric function proves the theorem. $\square$

\

\textbf{Fourth method.}

This method starts with the Rodrigues formula for the Jacobi polynomials
\[
P_n^{(\alpha,\beta)}(t)=\dfrac{(-1)^n}{2^n\, n!}\dfrac{1}{(1-t)^\alpha(1+t)^\beta}
\dfrac{d^n}{dt^n}\big((1-t)^\alpha(1+t)^\beta(1-t^2)^n \big).
\]
Application to the integral \eqref{2.000} gives
\begin{equation}
I=\int_{-1}^1\dfrac{(1-t)^\alpha(1+t)^\beta}{(z-t)^\lambda}P_n^{(\alpha,\beta)}(t)dt=
\dfrac{(-1)^n}{2^n\, n!}\int_{-1}^1\dfrac{1}{(z-t)^\lambda}\dfrac{d^n}{dt^n}
\big((1-t)^{\alpha+n}(1+t)^{\beta+n}\big)dt.
\label{2.7a}
 \end{equation}
One possibility is to apply partial integration. But we prefer to  apply the standard formula for the $n$-th order derivative of a product.
\[
\dfrac{d^n}{dt^n}\big(f(t)g(t)\big)=\sum_{k=0}^n\binom{n}{k}\dfrac{d^{n-k}}{dt^{n-k}}f(t)\dfrac{d^k}{dt^k}g(t).
\]
Then we get
\begin{align*}
&\dfrac{d^k}{dt^k}(1+t)^{\beta+n}=\dfrac{\Gamma(\beta+1+n)}{\Gamma(\beta+1+n-k)}(1+t)^{\beta+n-k} \\
&\dfrac{d^{n-k}}{dt^{n-k}}(1-t)^{\alpha+n}=(-1)^{n-k}\dfrac{\Gamma(\alpha+1+n)}{\Gamma(\alpha+1+k)}(1-t)^{\alpha+k}.
\end{align*}
Application to \eqref{2.7a} and interchanging the summation and the integral, which is allowed because the summation is convergent, gives
\[
I=\dfrac{(-1)^n}{2^n\, n!}\sum_{k=0}^n\binom{n}{k}(-1)^{n-k}
\dfrac{\Gamma(\alpha+1+n)}{\Gamma(\alpha+1+k)}\dfrac{\Gamma(\beta+1+n)}{\Gamma(\beta+1+n-k)}
\int_{-1}^1\dfrac{1}{(z-t)^\lambda}(1-t)^{\alpha+k}(1+t)^{\beta+n-k}dt.
\]
Using $t=2x-1$ changes the bounds from $(-1,1)$ to $(0,1)$. We use also $\binom{n}{k}=(-1)^k(-n)_k/k!$. We get 
\[
I=\dfrac{2^{\alpha+\beta+1}}{n!}\dfrac{1}{(z+1)^\lambda}\sum_{k=0}^n\dfrac{(-n)_k}{k!}
\dfrac{\Gamma(\alpha+1+n)}{\Gamma(\alpha+1+k)}\dfrac{\Gamma(\beta+1+n)}{\Gamma(\beta+1+n-k)}
\int_0^1 x^{\beta+n-k}(1-x)^{\alpha+k}\left(1-\dfrac{2}{z+1}x\right)^{-\lambda}dx.
\]
The integral gives a hypergeometric function
\begin{equation}
I=\dfrac{2^{\alpha+\beta+1}}{n!}\dfrac{1}{(z+1)^\lambda}
\dfrac{\Gamma(\alpha+1+n)\Gamma(\beta+1+n)}{\Gamma(\alpha+\beta+n+2)}
\sum_{k=0}^n\dfrac{(-n)_k}{k!}\hyp21{\lambda,1+\beta+n-k}{\alpha+\beta+n+2}{\dfrac{2}{z+1}}.
\label{2.8a}
\end{equation}
with $\Re(z) \geq 1$. Using a Gauss transformation results in
\[
I=\dfrac{2^{\alpha+\beta+1}}{n!}\dfrac{1}{(z-1)^\lambda}
\dfrac{\Gamma(\alpha+1+n)\Gamma(\beta+1+n)}{\Gamma(\alpha+\beta+n+2)}
\sum_{k=0}^n\dfrac{(-n)_k}{k!}\hyp21{\lambda,1+\alpha+k}{\alpha+\beta+n+2}{\dfrac{2}{1-z}}.
\]
Writing the hypergeometric function as a summation, using
\begin{equation}
(\alpha+1+k)_j=\dfrac{(\alpha+1+j)_k(\alpha+1)_j}{(\alpha+1)_k}
\label{2.8aa}
\end{equation}
and interchanging the summations, which is allowed because the summations are convergent, gives
\[
I=\dfrac{2^{\alpha+\beta+1}}{n!}\dfrac{1}{(z-1)^\lambda}
\dfrac{\Gamma(\alpha+1+n)\Gamma(\beta+1+n)}{\Gamma(\alpha+\beta+n+2)}
\sum_{j=0}^\infty\dfrac{(\lambda)_j(\alpha+1)_j}{(\alpha+\beta+n+2)_j}\dfrac{1}{j!}
\left(\dfrac{2}{1-z}\right)^j
\sum_{k=0}^n\dfrac{(-n)_k(\alpha+1+j)_k}{(\alpha+1)_k}\dfrac{1}{k!}.
\]
The last summation is standard. There rests
\begin{equation}
I=\dfrac{2^{\alpha+\beta+1}}{n!}\dfrac{1}{(z-1)^\lambda}
\dfrac{\Gamma(\alpha+1)\Gamma(\beta+1+n)}{\Gamma(\alpha+\beta+n+2)}
\sum_{j=n}^\infty\dfrac{(\lambda)_j(\alpha+1)_j}{(\alpha+\beta+n+2)_j}
\dfrac{\Gamma(n-j)}{\Gamma(-j)}\dfrac{1}{j!}\left(\dfrac{2}{1-z}\right)^j.
\label{2.8bb}
\end{equation}
Using the transformation  $k=j-n$ for the summation and \eqref{2.12aa} gives at last
\[
I=\dfrac{2^{\alpha+\beta+n+1}\Gamma(\beta+n+1)}{(z-1)^{n+\lambda}\Gamma(n+1)\Gamma(\lambda)}
\sum_{k=0}^\infty\dfrac{\Gamma(k+\alpha+n+1)\Gamma(k+n+\lambda)}{\Gamma(k+\alpha+\beta+2n+2)}\dfrac{1}{k!}
\left(\dfrac{2}{1-z}\right)^k.
\]
This is the same equation as \eqref{2.6a}. So this proves the theorem. $\square$

\

\textbf{Fifth method.}

This method starts with the first steps of the third method. However from equation \eqref{2.8a} we go into another direction. Writing the hypergeometric function as a summation and using
\[
(1+\beta+n-k)_j=(-1)^k(\beta+1+n)_{j-k}(-\beta-n)_k
\]
we get
\begin{multline*}
I=\int_{-1}^1\dfrac{(1-t)^\alpha(1+t)^\beta}{(z-t)^\lambda}P_n^{(\alpha,\beta)}(t)dt=
\dfrac{2^{\alpha+\beta+1}}{n!}\dfrac{1}{(z+1)^\lambda}
\dfrac{\Gamma(\alpha+1+n)\Gamma(\beta+1+n)}{\Gamma(\alpha+\beta+n+2)} \\
\sum_{j=0}^\infty\sum_{k=0}^n\dfrac{(\beta+1+n)_{j-k}(\lambda)_j(-n)_k(-\beta-n)_k}{(\alpha+\beta+n+2)_j}\dfrac{1}{j!k!}\left(\dfrac{2}{z+1}\right)^j(-1)^k
\end{multline*}
The double sum can be written as a Horn $H_2$ function \cite[5.7.1(14)]{5}
\begin{multline}
I=\int_{-1}^1\dfrac{(1-t)^\alpha(1+t)^\beta}{(z-t)^\lambda}P_n^{(\alpha,\beta)}(t)dt=
\dfrac{2^{\alpha+\beta+1}}{n!}\dfrac{1}{(z+1)^\lambda}
\dfrac{\Gamma(\alpha+1+n)\Gamma(\beta+1+n)}{\Gamma(\alpha+\beta+n+2)} \\
H_2\left(\beta+1+n,\lambda,-n,-\beta-n,\alpha+\beta+n+2;\dfrac{2}{z+1},-1\right)
\label{2.9}
\end{multline}
In \cite[(B5)]{4} we prove for $H_2(x,-1)$, defined as the limit of $H_2(x,y)$ for
$y \downarrow -1$, the following property
\begin{align*}
H_2(a_0,b_1,b_2,c_1,c_2;x,-1)&=\dfrac{\Gamma(1-a_0)\Gamma(1-a_0-b_2-c_1)}{\Gamma(1-a_0-b_2)\Gamma(1-a_0-c_1)}\hyp32{a_0+b_2,a_0+c_1,b_1}{a_0+b_2+c_1,c_2}{x}+ \\
&+\dfrac{\Gamma(1-a_0)\Gamma(a_0+b_2+c_1-1)}{\Gamma(b_2)\Gamma(c_1)}
\dfrac{\Gamma(c_2)\Gamma(b_1-a_0-b_2-c_1+1)}{\Gamma(b_1)\Gamma(c_2-a_0-b_2-c_1+1)} \\
&\qquad\qquad x^{1-a_0-b_2-c_1}
\hyp32{1-b_2,1-c_1,b_1-a_0-b_2-c_1+1}{2-a_0-b_2-c_1,c_2-a_0-b_2-c_1+1}{x}
\end{align*}
with the condition $\Re(b_1-a_0-b_2-c_1+1)>0$ and $0 \leq x<1$. Application of
\[
a_0=\beta+1+n, \qquad b_1=\lambda, \qquad b_2=-n, \qquad c_1=-\beta-n-\epsilon, \qquad c_2=\alpha+\beta+n+2
\]
with $\epsilon \rightarrow 0$ gives
\begin{align*}
H_2&\left(\beta+1+n,\lambda,-n,-\beta-n,\alpha+\beta+n+2;\dfrac{2}{z+1},-1\right)= \\
&=\lim_{\epsilon \rightarrow 0}\dfrac{\Gamma(-\beta-n)}{\Gamma(-\beta)}
\dfrac{\Gamma(n+\epsilon)\Gamma(1-\epsilon-n)}{\Gamma(\epsilon)}\dfrac{1}{\Gamma(1-\epsilon-n)}
\hyp32{\beta+1,1,\lambda}{1-\epsilon-n,\alpha+\beta+n+2,}{\dfrac{2}{z+1}}+ \\
&+\dfrac{\Gamma(\alpha+\beta+n+2)\Gamma(n+\lambda)}{\Gamma(\lambda)\Gamma(\alpha+\beta+2n+2)}\left(\dfrac{2}{z+1}\right)^n
\hyp21{\beta+n+1,n+\lambda}{\alpha+\beta+2n+2}{\dfrac{2}{z+1}}
\end{align*}
with $z>1$. For the quotient of the Gamma functions we have
\[
\dfrac{\Gamma(-\beta-n)}{\Gamma(-\beta)}=(-1)^n\dfrac{1}{(\beta+1)_n} \qquad\qquad
\lim_{\epsilon \rightarrow 0}\dfrac{\Gamma(n+\epsilon)\Gamma(1-\epsilon-n)}{\Gamma(\epsilon)}=(-1)^n
\]
Applying \eqref{2.002} with $M=n+\epsilon-1=n-1$ we get after some manipulations with the Gamma functions and the Pochhammer symbols
\begin{align*}
H_2&\left(\beta+1+n,\lambda,-n,-\beta-n,\alpha+\beta+n+2;\dfrac{2}{z+1},-1\right)= \\
&=\dfrac{1}{(\beta+1)_n}\left(\dfrac{2}{z+1}\right)^n
\dfrac{(\beta+1)_n(1)_n(\lambda)_n}{\Gamma(n+1)(\alpha+\beta+n+2)_n}
\hyp21{\beta+1+n,\lambda+n}{\alpha+\beta+2n+2}{\dfrac{2}{z+1}}+ \\
&+\dfrac{(\lambda)_n}{(\alpha+\beta+n+2)_n}\left(\dfrac{2}{z+1}\right)^n
\hyp21{\beta+n+1,n+\lambda}{\alpha+\beta+2n+2}{\dfrac{2}{z+1}}
\end{align*}
Simplification gives
\begin{multline*}
H_2\left(\beta+1+n,\lambda,-n,-\beta-n,\alpha+\beta+n+2;\dfrac{2}{z+1},-1\right)= \\
=\dfrac{2^{n+1}(\lambda)_n}{(\alpha+\beta+n+2)_n}\left(\dfrac{1}{z+1}\right)^n
\hyp21{\beta+n+1,n+\lambda}{\alpha+\beta+2n+2}{\dfrac{2}{z+1}}
\end{multline*}
Application to \eqref{2.9} results in
\begin{multline*}
\int_{-1}^1\dfrac{(1-t)^\alpha(1+t)^\beta}{(z-t)^\lambda}P_n^{(\alpha,\beta)}(t)dt= \\
=\dfrac{2^{\alpha+\beta+n+2}(\lambda)_n}{n!}B(\alpha+1+n,\beta+1+n)\dfrac{1}{(z+1)^{n+\lambda}}
\hyp21{\beta+n+1,n+\lambda}{\alpha+\beta+2n+2}{\dfrac{2}{z+1}}
\end{multline*}
Applying a standard Gauss transformation proves the theorem. $\square$

\

\textbf{Sixth method.}

We call this method a brute force method. The Jacobi polynomial is written as a hypergeometric function which is written as a summation. After interchanging the summation and the integration which is allowed because of the convergence of the Jacobi polynomial and integration gives a double summation which can be simplified until we reach the final result.
\begin{align*}
I&=\int_{-1}^1\dfrac{(1-t)^\alpha(1+t)^\beta}{(z-t)^\lambda}P_n^{(\alpha,\beta)}(t)dt  \\
&=\int_{-1}^1\dfrac{(1-t)^\alpha(1+t)^\beta}{(z-t)^\lambda}\binom{n+\alpha}{n}
\hyp21{-n,n+\alpha+\beta+1}{\alpha+1}{\dfrac{1-t}{2}} \\
&=\dfrac{\Gamma(n+\alpha+1)}{\Gamma(n+1)\Gamma(\alpha_n+1)}
\sum_{i=0}^n\dfrac{(-n)_i(n+\alpha+\beta+1)_i}{(\alpha+1)_i}\dfrac{1}{i!}
\left(\dfrac{1}{2}\right)^i
\int_{-1}^1\dfrac{(1-t)^{\alpha+i}(1+t)^\beta}{(z-t)^\lambda}dt
\end{align*}
Integration gives a hypergeometric function which can be written as a summation
\[
I=2^{\alpha+\beta+1}\dfrac{\Gamma(n+\alpha+1)\Gamma(\beta+1)}
{\Gamma(n+1)\Gamma(\alpha+\beta+2)}\left(\dfrac{1}{z-1}\right)^\lambda
\sum_{i=0}^n\dfrac{(-n)_i(n+\alpha+\beta+1)_i}{(\alpha+\beta+2)_i}\dfrac{1}{i!}
\sum_{j=0}^\infty\dfrac{(\lambda)_j(\alpha+1+i)_j}{(\alpha+\beta+2+i)_j}\dfrac{1}{j!}
\left(\dfrac{2}{1-z}\right)^j
\]
Making use of \eqref{2.8aa} gives
\begin{multline*}
I=2^{\alpha+\beta+1}\dfrac{\Gamma(n+\alpha+1)\Gamma(\beta+1)}
{\Gamma(n+1)\Gamma(\alpha+\beta+2)}\left(\dfrac{1}{z-1}\right)^\lambda
\sum_{j=0}^\infty\dfrac{(\alpha+1)_j(\lambda)_j}{(\alpha+\beta+2)_j}\
\dfrac{1}{j!}\left(\dfrac{2}{1-z}\right)^j \\
\sum_{i=0}^n\dfrac{(-n)_i(n+\alpha+\beta+1)_i(\alpha+1+j)_i}
{(\alpha+1)_i(\alpha+\beta+2+j)_i}\dfrac{1}{i!}
\end{multline*}
For the last summation we can use Saalsch\"utz's theorem. After much manipulations with the Gamma functions we get
\[
I=2^{\alpha+\beta+1}\dfrac{\Gamma(\alpha+1)\Gamma(\beta+1+n)}
{\Gamma(n+1)\Gamma(\alpha+\beta+2+n)}\left(\dfrac{1}{z-1}\right)^\lambda
\sum_{j=0}^\infty\dfrac{(\alpha+1)_j(\lambda)_j}{(\alpha+\beta+2+n)_j}
\dfrac{1}{j!}\left(\dfrac{2}{1-z}\right)^j \dfrac{\Gamma(n-j)}{\Gamma(-j)}
\]
This is the same equation as \eqref{2.8bb}, so this proves the theorem. $\square$

\

\

\textbf{\fontsize{10.5}{12.5}\selectfont Proof of Theorem 2.}

To prove the integral we use the brute-force method. Set
\[
I=\int_x^1\dfrac{(1-t)^\alpha(1+t)^\beta}{(t-x)^\lambda}P_n^{(\alpha,\beta)}(t)dt
\]
 We write the Jacobi polynomial as a hypergeometric function. Writing this hypergeometric function as a summation and interchanging the integral and the summation (which is allowed because of the convergence of the integral and the summation) results in
\[
I=\dfrac{\Gamma(n+\alpha+1)}{\Gamma(n+1)\Gamma(\alpha+1)}
\sum_{i=0}^n\dfrac{(-n)_i(n+\alpha+\beta+1)_i}{(\alpha+1)_i}\dfrac{1}{i!}
\left(\dfrac{1}{2}\right)^i
\int_x^1\dfrac{(1-t)^{\alpha+i}(1+t)^\beta}{(t-x)^\lambda}dt
\]
The indefinite integral can be computed and there arises an $F_1$ Appell function  \cite[5.7.1(6)]{5}.
\begin{multline*}
I=\dfrac{\Gamma(n+\alpha+1)}{\Gamma(n+1)\Gamma(\alpha+1)}
\sum_{i=0}^n\dfrac{(-n)_i(n+\alpha+\beta+1)_i}{(\alpha+1)_i}\dfrac{1}{i!}
\left(\dfrac{1}{2}\right)^i \\
\left[-\dfrac{2^\beta}{(\alpha+1+i)}\dfrac{(1-t)^{\alpha+1+i}}{(1-x)^\lambda}
F_1\left(\begin{array}{l}
	\alpha+1+i,-\beta,\lambda \\
	\alpha+2+i
\end{array};\dfrac{1-t}{2},\dfrac{1-t}{1-x}\right)
\right]^1_x
\end{multline*}
Substitution of the boundary conditions gives after some simplification
\begin{multline*}
I=-2^\beta\dfrac{\Gamma(n+\alpha+1)}{\Gamma(n+1)\Gamma(\alpha+2)}(1-x)^{\alpha+1-\lambda} \\
\sum_{i=0}^n\dfrac{(-n)_i(n+\alpha+\beta+1)_i}{(\alpha+2)_i}\dfrac{1}{i!}
\left(\dfrac{1-x}{2}\right)^i 
F_1\left(\begin{array}{l}
	\alpha+1+i,-\beta,\lambda \\
	\alpha+2+i
\end{array};\dfrac{1-x}{2},1\right)
\end{multline*}
For the Appell function with one of the arguments equal to $1$ we have \cite[5.10(10)]{5}
\begin{equation}
F_1\left(\begin{array}{l}
a,b_1,b_2 \\
c
\end{array};x,1\right)
=\dfrac{\Gamma(c)\Gamma(c-a-b_2)}{\Gamma(c-a)\Gamma(c-b_2)}
\hyp32{a,b_1,a+1-c}{c,a+b_2-c+1}{x}
\label{2.15a}
\end{equation}
Interchanging $(b_1,b_2)$ into $(b_2,b_1)$ and $(x,1)$ into $(1,x)$ results in
\begin{equation}
F_1\left(\begin{array}{l}
	a,b_1,b_2 \\
	c
\end{array};1,x\right)
=\dfrac{\Gamma(c)\Gamma(c-a-b_1)}{\Gamma(c-a)\Gamma(c-b_1)}
\hyp32{a,b_2,a+1-c}{c,a+b_1-c+1}{x}
\label{2.15b}
\end{equation}

Substitution gives
\begin{multline}
I=-2^\beta\dfrac{\Gamma(n+\alpha+1)}{\Gamma(n+1)\Gamma(\alpha+2)}(1-x)^{\alpha+1-\lambda} \\
\sum_{i=0}^n\dfrac{(-n)_i(n+\alpha+\beta+1)_i}{(\alpha+2)_i}\dfrac{1}{i!}
\left(\dfrac{1-x}{2}\right)^i 
\dfrac{\Gamma(\alpha+2+i)\Gamma(1-\lambda)}{\Gamma(\alpha+2-\lambda+i)}
\hyp21{\alpha+1+i,-\beta}{\alpha+2-\lambda+i}{\dfrac{1-x}{2}}
\label{2.13}
\end{multline}
Using a standard transformation for the hypergeometric function and writing it as a summation gives after some simplification 
\begin{multline*}
I=-2^\beta\dfrac{\Gamma(n+\alpha+1)\Gamma(1-\lambda)}{\Gamma(n+1)\Gamma(\alpha-\lambda+2)}(1+x)^\beta (1-x)^{\alpha+1-\lambda} \\
\sum_{i=0}^n\dfrac{(-n)_i(n+\alpha+\beta+1)_i}{(\alpha-\lambda+2)_i}\dfrac{1}{i!}
\left(\dfrac{1-x}{2}\right)^i 
\sum_{j=0}^\infty\dfrac{(1-\lambda)_j(-\beta)_j}{(\alpha+2-\lambda+i)_j}\dfrac{1}{j!}
\left(\dfrac{x-1}{x+1}\right)^j
\end{multline*}
Making use of $(\alpha-\lambda+2)_i(\alpha+2-\lambda+i)_j=(\alpha-\lambda+2)_{i+j}$ gives
\begin{multline*}
I=-2^\beta\dfrac{\Gamma(n+\alpha+1)\Gamma(1-\lambda)}{\Gamma(n+1)\Gamma(\alpha-\lambda+2)}(1+x)^\beta (1-x)^{\alpha+1-\lambda} \\
\sum_{i=0}^n\sum_{j=0}^\infty
\dfrac{(-n)_i(-\beta)_j(n+\alpha+\beta+1)_i(1-\lambda)_j}{(\alpha-\lambda+2)_{i+j}}
\dfrac{1}{i!j!}\left(\dfrac{1-x}{2}\right)^i \left(\dfrac{x-1}{x+1}\right)^j
\end{multline*}	
The double summation is an $F_3$ \text{\em zero-balanced} Appell function.  We call the $F_3$ functions zero-balanced because in this case the sum of the values of the upper parameters is equal to the sum of the lower parameters. $\cite[5.7.1(8)]{5}$.
\[
I=-2^\beta\dfrac{\Gamma(n+\alpha+1)\Gamma(1-\lambda)}{\Gamma(n+1)\Gamma(\alpha-\lambda+2)}(1+x)^\beta (1-x)^{\alpha+1-\lambda}
F_3\left(\begin{array}{l}
	-n,-\beta,\alpha+\beta+1+n,1-\lambda \\
	\alpha-\lambda+2
\end{array};\dfrac{1-x}{2},\dfrac{x-1}{x+1}\right)
\]
The $F_3$ Appell function can be converted into a hypergeometric function. \cite[p. 302 (89)]{9} gives 
\begin{equation}
F_3\left(\begin{array}{l}
	a,b,c,d \\
	a+b+c+d
\end{array};\dfrac{1-x}{2},\dfrac{x-1}{x+1}\right)
=\left(\dfrac{1+x}{2}\right)^b \hyp21{a+b,b+c}{a+b+c+d}{\dfrac{1-x}{2}}
\label{3.17}
\end{equation}
Application results at last in
\[
\int_x^1\dfrac{(1-t)^\alpha(1+t)^\beta}{(t-x)^\lambda}P_n^{(\alpha,\beta)}(t)dt 
=2^\beta\dfrac{\Gamma(n+\alpha+1)\Gamma(1-\lambda)}{\Gamma(n+1)\Gamma(\alpha-\lambda+2)}(1-x)^{\alpha+1-\lambda}\hyp21{\alpha+n+1,-\beta-n}{\alpha-\lambda+2}{\dfrac{1-x}{2}}
\]
This proves the theorem. $\square$

\

\

\textbf{Remark.}\quad
Taking the limit for $x \downarrow -1$ on both sides of \eqref{2.13} gives
\begin{align*}
&\int_{-1}^1(1-t)^\alpha(1+t)^{\beta-\lambda}P_n^{(\alpha,\beta)}(t)dt= \\
&=2^{\alpha+\beta-\lambda+1}
\dfrac{\Gamma(n+\alpha+1)}{\Gamma(n+1)\Gamma(\alpha+2)}
\sum_{i=0}^n\dfrac{(-n)_i(n+\alpha+\beta+1)_i}{(\alpha+2)_i}\dfrac{1}{i!}
\dfrac{\Gamma(\alpha+2+i)\Gamma(1-\lambda)}{\Gamma(\alpha+2-\lambda+i)}
\hyp21{\alpha+1+i,-\beta}{\alpha+2-\lambda+i}{1} \\
&=2^{\alpha+\beta-\lambda+1}
\dfrac{\Gamma(n+\alpha+1)\Gamma(\beta-\lambda+1)}{\Gamma(n+1)\Gamma(\alpha+\beta+2-\lambda)}
\sum_{i=0}^n\dfrac{(-n)_i(n+\alpha+\beta+1)_i}{(\alpha+\beta+2-\lambda)_i}\dfrac{1}{i!} \\
&=2^{\alpha+\beta-\lambda+1}\binom{-\lambda}{n}B(\alpha+n+1,\beta-\lambda+1)
\end{align*}
with conditions:  $\alpha>-1,\beta-\lambda>-1$. This integral is known \cite[7.39.3]{2}.

\

\

\textbf{\fontsize{10.5}{12.5}\selectfont Proof of Theorem 3.}

No prove is needed. Looking at the integral the Jacobi polynomial can be converted into a Gegenbauer polynomial. If we substitute for $\alpha$ and $\beta$ the value $\alpha-1/2$ and using the properties
\begin{align*}
&P_n^{(\alpha-1/2,\alpha-1/2)}(t)=\dfrac{\Gamma(2\alpha)\Gamma\left(\alpha+n+\dfrac{1}{2}\right)}{\Gamma(2\alpha+n)\Gamma\left(\alpha+\dfrac{1}{2}\right)}C_n^{(\alpha)}(t)
\qquad\qquad\qquad\qquad\ \ \cite[18.7.1]{11} \\
&(-1)^n\dfrac{\Gamma(\lambda+n)\Gamma(1-\lambda-n)}{\Gamma(\lambda)}=\Gamma(1-\lambda) \\
&\hyp21{a,1-a}{b}{\dfrac{1}{2}}=2^{1-b}\sqrt{\pi}
\dfrac{\Gamma(b)}{\Gamma\left(\dfrac{a+b}{2}\right)\Gamma\left(\dfrac{1+b-a}{2}\right)}
\qquad\qquad\qquad \cite[7.3.7. (8)]{10}
\end{align*}
we arrive at equation \eqref{4}.

\

\

\textbf{\fontsize{10.5}{12.5}\selectfont Proof of Theorem 4.}

To prove this theorem we write the Jacobi polynomial as a hypergeometric function.
\[
I=\int_{-1}^x\dfrac{(1-t)^\alpha(1+t)^\beta}{(t-x)^\lambda}P_n^{(\alpha,\beta)}(t)dt=
\int_{-1}^x\dfrac{(1-t)^\alpha(1+t)^\beta}{(t-x)^\lambda}\binom{n+\alpha}{n}
\hyp21{-n,n+\alpha+\beta+1}{\alpha+1}{\dfrac{1-t}{2}}dt
\]
Writing the hypergeometric function as a summation and interchanging the integral and the summation (which is allowed because of the convergence of the integral and the summation) results in
\[
I=\dfrac{\Gamma(n+\alpha+1)}{\Gamma(n+1)\Gamma(\alpha+1)}
\sum_{i=0}^n\dfrac{(-n)_i(n+\alpha+\beta+1)_i}{(\alpha+1)_i}\dfrac{1}{i!}
\left(\dfrac{1}{2}\right)^i
\int_{-1}^x\dfrac{(1-t)^\alpha(1+t)^\beta}{(t-x)^\lambda}dt
\]
The integral is known (see the proof of Theorem 2). So we get
\begin{multline*}
I=2^\beta\dfrac{\Gamma(n+\alpha+1)}{\Gamma(n+1)\Gamma(\alpha+2)}(1-x)^{-\lambda}
\sum_{i=0}^n\dfrac{(-n)_i(n+\alpha+\beta+1)_i}{(\alpha+2)_i}\dfrac{1}{i!}
\left(\dfrac{1}{2}\right)^i \\
\left[2^{\alpha+1+i}
F_1\left(\begin{array}{l}
	\alpha+1+i,-\beta,\lambda \\
	\alpha+2+i
\end{array};1,\dfrac{2}{1-x}\right)
-(1-x)^{\alpha+1+i}
F_1\left(\begin{array}{l}
	\alpha+1+i,-\beta,\lambda \\
	\alpha+2+k
\end{array};\dfrac{1-x}{2},1\right)
\right]
\end{multline*}
Using \eqref{2.15a} and \eqref{2.15b} gives after a lot of manipulations with the Gamma functions and the Pochhammer symbols
\begin{align*}
I&=\dfrac{2^{\alpha+\beta+1}}{(1-x)^\lambda}\left(\dfrac{x+1}{x-1}\right)^{-\lambda}
\dfrac{\Gamma(n+\alpha+1)\Gamma(\beta+1)}{\Gamma(n+1)\Gamma(\alpha+\beta+2)}
\sum_{j=0}^\infty\dfrac{(\beta+1)_j(\lambda)_j}{(\alpha+\beta+2)_j}\dfrac{1}{j!}
\left(\dfrac{2}{x+1}\right)^j
\sum_{i=0}^n\dfrac{(-n)_i(n+\alpha+\beta+1)_i}{(\alpha+\beta+2+j)_i}\dfrac{1}{i!}- \\
&-\dfrac{(1-x)^{\alpha+1}}{(1-x)^\lambda}(x+1)^\beta
\dfrac{\Gamma(n+\alpha+1)\Gamma(1-\lambda)}{\Gamma(n+1)\Gamma(\alpha-\lambda+2)}
\sum_{i=0}^n\sum_{j=0}^\infty
\dfrac{(-n)_i(-\beta)_j(n+\alpha+\beta+1)_i(1-\lambda)_i}{(\alpha-\lambda+2)_{i+j}}
\dfrac{1}{i!j!}\left(\dfrac{1-x}{2}\right)^i\left(\dfrac{x-1}{x+1}\right)^j
\end{align*}
The last summation in the first term can be evaluated. The double summation in the second term is an $F_3$ Appell function.
\begin{align*}
I&=\dfrac{2^{\alpha+\beta+1}}{(1-x)^\lambda}\left(\dfrac{x+1}{x-1}\right)^{-\lambda}
\dfrac{\Gamma(n+\alpha+1)\Gamma(\beta+1)}{\Gamma(n+1)\Gamma(\alpha+\beta+2+n)}
\sum_{j=n}^\infty\dfrac{(\beta+1)_j(\lambda)_j}{(\alpha+\beta+2+n)_j}\dfrac{1}{j!}
\left(\dfrac{2}{x+1}\right)^j\dfrac{j+1}{j-n+1}- \\
&-\dfrac{(1-x)^{\alpha+1}}{(1-x)^\lambda}(x+1)^\beta
\dfrac{\Gamma(n+\alpha+1)\Gamma(1-\lambda)}{\Gamma(n+1)\Gamma(\alpha-\lambda+2)}
F_3\left(\begin{array}{l}
	-n,-\beta,n+\alpha+\beta+1,1-\lambda \\
	\alpha-\lambda+2
\end{array};\dfrac{1-x}{2},\dfrac{x-1}{x+1}\right)
\end{align*}
Setting $j=m+n$ and $(a)_{i+j}=(a+i)_j(a)_i$ for the summation in the first term we get
\begin{align*}
I&=\dfrac{2^{\alpha+\beta+1}(\lambda)_n}{(1-x)^\lambda}\left(\dfrac{x+1}{x-1}\right)^{-\lambda}
\left(\dfrac{2}{x+1}\right)^n
\dfrac{\Gamma(n+\alpha+1)\Gamma(\beta+n+1)}{\Gamma(n+1)\Gamma(\alpha+\beta+2+2n)}
\sum_{m=0}^\infty\dfrac{(\beta+1+n)_m(\lambda+n)_m}{(\alpha+\beta+2+2n)_m}\dfrac{1}{m!}
\left(\dfrac{2}{x+1}\right)^m- \\
&-(1-x)^{\alpha-\lambda+1}(x+1)^\beta
\dfrac{\Gamma(n+\alpha+1)\Gamma(1-\lambda)}{\Gamma(n+1)\Gamma(\alpha-\lambda+2)}
F_3\left(\begin{array}{l}
	-n,-\beta,n+\alpha+\beta+1,1-\lambda \\
	\alpha-\lambda+2
\end{array};\dfrac{1-x}{2},\dfrac{x-1}{x+1}\right)
\end{align*}
Writing the summation in the first term as a hypergeometric function we get after some simplification
\begin{align*}
I&=\dfrac{2^{\alpha+\beta+n+1}}{(1-x)^{n+\lambda}}\binom{-\lambda}{n}
B(\alpha+n+1,\beta+n+1)\hyp21{\alpha+n+1,n+\lambda}{\alpha+\beta+2+2n}{\dfrac{2}{1-x}}- \\
&-(1-x)^{\alpha-\lambda+1}(x+1)^\beta
\dfrac{\Gamma(n+\alpha+1)\Gamma(1-\lambda)}{\Gamma(n+1)\Gamma(\alpha-\lambda+2)}
F_3\left(\begin{array}{l}
	-n,-\beta,n+\alpha+\beta+1,1-\lambda \\
	\alpha-\lambda+2
\end{array};\dfrac{1-x}{2},\dfrac{x-1}{x+1}\right)
\end{align*}
For the $F_3$ Appell function we use \eqref{3.17} and get
\begin{align*}
I&=\dfrac{2^{\alpha+\beta+n+1}}{(1-x)^{n+\lambda}}\binom{-\lambda}{n}
B(\alpha+n+1,\beta+n+1)\hyp21{\alpha+n+1,n+\lambda}{\alpha+\beta+2+2n}{\dfrac{2}{1-x}}- \\
&-2^\beta\dfrac{\Gamma(n+\alpha+1)\Gamma(1-\lambda)}{\Gamma(n+1)\Gamma(\alpha-\lambda+2)}
(1-x)^{\alpha+1-\lambda}\hyp21{\alpha+n+1,-\beta-n}{\alpha-\lambda+2 }{\dfrac{1-x}{2}}
\end{align*}
This proves the theorem. $\square$

\


\begin{thebibliography}{7}
%
\bibitem{4}
Diekema, E. {\em The fractional orthogonal derivative for functions of one and two variables.} PhD Thesis,\qquad\qquad  University of Amsterdam, (2018), 
https://hdl.handle.net/11245.1/a6ed8a3f-0831-4f9a-9476-2661ec7e1f92.
%
\bibitem{6}
Diekema, E. {\em A correlation function for the classical orthogonal polynomials.}
arXiv:2011.07498 (2020).
%
\bibitem{5}
Erd\'elyi, A. {\em Higher transcendental functions, Vol. I}. McGraw-Hill, 1953.
%
\bibitem{7}
Erd\'elyi, A. {\em Table of integral transforms, Vol. II}. McGraw-Hill, 1954.
%
\bibitem{9}
Srivastava, H.M, Karlsson, P.W.  {\em Multiple Gaussian hypergeometric series}.
Ellis Horwood Limited, 1985.
%
\bibitem{2}
Gradshteyn, I.S. Ryzhik, I.M. {\em Table of Integrals, Series and Products}. Eighth ed. Elsevier 2014.
%
\bibitem{11}
Olver, F.W.J. et al. {\em NIST Handbook of mathematical functions}.
Cambridge University Press, 2010.  http://dmlf.nist.gov 
%
\bibitem{3}
Prudnikov, A.P. Brychkov, Y.A. Marichev, O.I. {\em Integrals and Series} Vol.2. Gordon and Breach Science Publishers 1986.
\bibitem{10}
Prudnikov, A.P. Brychkov, Y.A. Marichev, O.I. {\em Integrals and Series} Vol.3. Gordon and Breach Science Publishers 1990.
%
\bibitem{1}
Szmytkowski, R. {\em Some integrals and series involving the Gegenbauer polynomials and the Legendre functions on the cut (-1,1)}. arXiv:1107.2680v2 (2011).
\end{thebibliography}
\end{document}